\documentclass[12pt]{article}
\usepackage{amssymb,epsfig,times,graphicx}
\usepackage{amsmath}
\usepackage{graphicx}

\newtheorem{thm}{Theorem}

\newtheorem{prop}[thm]{Proposition}

\newtheorem{defi}[thm]{Definition}

\newtheorem{cor}[thm]{Corollary}

\newtheorem{rmk}[thm]{Remark}

\newcommand{\beq}{\begin{equation}}
\newcommand{\eeq}{\end{equation}}
\newcommand{\de}{\partial}
\def\d{\partial}
\def\n{\noindent}

\newcommand{\lm}{\lambda}

\begin{document}
\title{$F$-manifolds and integrable systems\\ of hydrodynamic type}
\author{Paolo Lorenzoni${}^{*}$, Marco Pedroni${}^{**}$, Andrea Raimondo${}^{***}$\\
\\
{\small * Dipartimento di Matematica e Applicazioni}\\
{\small Universit\`a di Milano-Bicocca}\\
{\small Via Roberto Cozzi 53, I-20125 Milano, Italy}\\
{\small paolo.lorenzoni@unimib.it}\\
\\
{\small ** Dipartimento di Ingegneria dell'informazione e metodi matematici}\\
{\small Universit\`a di Bergamo - Sede di Dalmine}\\
{\small viale Marconi 5, I-24044 Dalmine BG, Italy}\\
{\small marco.pedroni@unibg.it}\\
\\
{\small *** Department of Mathematics, Imperial College}\\
{\small 180 Queen's Gate, London SW7 2AZ, UK}\\
{\small a.raimondo@imperial.ac.uk}}
\date{}

\maketitle

\begin{abstract}
We investigate the role of Hertling-Manin condition on the structure constants of an associative commutative algebra
 in the theory of integrable systems of hydrodynamic type. In such a framework
  we introduce the notion of $F$-manifold with compatible connection generalizing
   a structure introduced by Manin.
\end{abstract}

\section{Introduction}
In their seminal papers \cite{dn89,ts91}, Dubrovin, Novikov, and Tsarev pointed out a deep relation between the integrability properties of systems of PDEs of hydrodynamic type
\beq\label{hdf}
u^i_{t}=V^i_j u^j_x,\qquad i=1,\dots,n,
\eeq
(sum over repeated indices is understood) and geometrical---in particular,
Riemannian---structures on the target manifold $M$, where $(u^1,\dots,u^n)$ play the role of coordinates. Probably, the most important of such structures is the notion of Frobenius manifold,
 introduced by Dubrovin (see, e.g., \cite{du93}) in order to give a coordinate-free description of the famous WDVV equations. A crucial
  ingredient involved in the definition of Frobenius manifolds is a $(1,2)$-type tensor field $c$ giving an associative commutative product on every tangent space:
$$
(X\circ Y)^i:=c^i_{jk}X^j Y^k\ ,
$$
where $X$ and $Y$ are vector fields. More recently \cite{HM}, Hertling and Manin showed that this product satisfies the condition
\beq
\label{hmcond}
\begin{aligned}
&&[X\circ Y,Z\circ W]-[X\circ Y,Z]\circ W-[X\circ Y,W]\circ Z-X\circ[Y,Z\circ W]+X\circ[Y,Z]\circ W+\\
&&+X\circ[Y,W]\circ Z-Y\circ[X,Z\circ W]+Y\circ[X,Z]\circ W+Y\circ[X,W]\circ Z=0\ ,
\end{aligned}
\eeq
or, in terms of the components of $c$,
\beq
\label{hmcomponents}
(\d_s c^k_{jl})c^s_{im}+(\d_j c^s_{im})c^k_{sl}-(\d_s c^k_{im})c^s_{jl}-(\d_i c^s_{jl})c^k_{sm}-(\d_l c^s_{jm})c^k_{si}-(\d_m c^s_{li})c^k_{js}=0\ .
\eeq
They called \emph{$F$-manifold} a manifold endowed with an associative commutative multiplicative structure satisfying condition \eqref{hmcond}.

The aim of this paper is to study the properties of the PDEs of hydrodynamic type
 associated with $F$-manifolds. The system (\ref{hmcomponents}) and its relation with integrable systems has been considered from a different point of view in \cite{KM}. Here, following the insights coming from the case of the principal hierarchy in the context of
Frobenius manifolds, we will assume such PDEs to be of the form
\beq\label{HTS}
u^i_{t}=(V_{X})^i_j u^j_x,\,\,\,i=1,\dots,n\ ,\qquad (V_{X})^i_j:=c^i_{jk}X^k,
\eeq
where $X$ is a vector field on $M$ and $c$ satisfies (\ref{hmcond}). These assumptions have two important consequences, spelled out respectively in Section 2 and 3:

1. For any choice of the vector field $X$, the Haantjes tensor associated with the (1,1) tensor field $V_{X}$ vanishes.
 
2. They allow one to write the condition of commutativity of two flows of the form (\ref{HTS}) as a simple requirement on the corresponding vector fields on $M$.

Starting from Section 4, we put into the game an additional structure, namely a connection $\nabla$ satisfying the symmetry condition
\begin{equation}\label{sccintro}
\left(\nabla_X c\right)\left(Y,Z\right)=\left(\nabla_Y c\right)\left(X,Z\right),
\end{equation}
for all vector fields $X$, $Y$, and $Z$. Remarkably, as shown by Hertling \cite{hert}, condition (\ref{hmcond}) follows from (\ref{sccintro}).

In Section 4, following Manin \cite{manin}, we study the special case where the connection $\nabla$ is flat and 
 we show how to construct an integrable hierarchy of hydrodynamic type. The costruction is divided in two steps. First---using a basis of flat vector fields---one  defines a set of flows, known as \emph{primary flows}. Then, from these flows one can define recursively the higher flows of the hierarchy. In this way, each primary flow turns out to be the starting point of a hierarchy. This construction is a straightforward generalization of the principal hierarchy defined by Dubrovin in the case of Frobenius manifolds \cite{du93}.

The general (non-flat) case is studied in Section 5, where we introduce the notion of $F$-manifold with compatible
 (non-flat) connection $\nabla$ and we show that the associated integrable systems of hydrodynamic type are defined by a family of vector fields satisfying the following condition:
\beq\label{admvfintro}
c^i_{jm}\nabla_k X^m=c^i_{km}\nabla_j X^m.
\eeq
In the non-flat case the existence of solutions of the above system is not guaranteed. Indeed, we prove that every solution $X$ of (\ref{admvfintro}) satisfies the condition
$$(R^k_{lmi}c^n_{pk}+R^k_{lip}c^n_{mk}+R^k_{lpm}c^n_{ik})X^l=0,$$
where $R$ is the curvature tensor of $\nabla$. It is thus natural to introduce the following requirement on the curvature:
\beq\label{shcintro}
R^k_{lmi}c^n_{pk}+R^k_{lip}c^n_{mk}+R^k_{lpm}c^n_{ik}=0\ .
\eeq
If the structure constants $c^i_{jk}$ admit canonical coordinates,  condition (\ref{shcintro}) is related to the well-known semi-Hamiltonian property introduced by Tsarev  \cite{ts91} as compatibilty condition for the linear system providing the symmetries
 of a diagonal system of hydrodynamic type.

In Section 6, motivated by the Hamiltonian theory of systems of hydrodynamic type, we consider the case of metric connections and we introduce
 the notion of Riemannian $F$-manifold. Finally, in Section 7, we discuss in details an important example: the reductions of the dispersionless KP hierarchy (also known as Benney chain).

\section{The Haantjes tensor}
\label{sec:haantjes}
An important class of systems of hydrodynamic type, widely studied in the literature, consists in those systems which admit diagonal form.  We say that a system \eqref{hdf} is diagonalizable if there exists a set of coordinates $\left(r^1,\dots,r^n\right)$---usually called \emph{Riemann invariants}---such that the tensor $V^i_j$ is diagonal in these coordinates: $V^i_j(r)=v^i\,\delta^i_j$. Then the system takes the (diagonal) form
\begin{equation*}
r^i_{t}=v^i(r^1,\dots,r^n) r^i_x,\,\,\,i=1,\dots,n\ .
\end{equation*}
It is important to recall that there exists an invariant criterion for the diagonalizability. One first introduces the \emph{Nijenhuis tensor} of $V$ as
$$N_V(X,Y)=[VX,VY]-V\,[X,VY]-V\,[VX,Y]+V^2\,[X,Y],$$
where $X$ and $Y$ are arbitrary vector fields, and then defines the \emph{Haantjes tensor} as
$$H_V(X,Y)=N_V(VX,VY)-V\,N_V(X,VY)-V\,N_V(VX,Y)+V^2N_V(X,Y).$$
In the case when $V$ has mutually distinct eigenvalues, then $V$ is diagonalizable if and only if its Haantjes tensor is identically zero. In this section, we consider the Haantjes tensor of
\begin{equation}\label{vcx}
(V_Z)^i_j=c^i_{jk}Z^k,
\end{equation}
where $c$ satisfies the Hertling-Manin condition \eqref{hmcond}. For a $(1,1)$- type tensor field of the form \eqref{vcx}, the Nijenhuis tensor reads
$$N_{V_Z}(X,Y)=[Z\circ X,Z\circ Y]+Z^2\circ [X,Y]-Z\circ [X,Z\circ Y]-Z\circ [Z\circ X,Y].$$
By using the Hertling-Manin condition (\ref{hmcond}) evaluated at $X=Z$, this can be written as
$$N_{V_Z}(X,Y)=[X\circ Z,Z]\circ Y-[X,Y]\circ Z\circ Y+[Z,Y\circ Z]\circ X-[Z,Y]\circ X\circ Z,$$
using this identity it is easy to prove the following
\begin{thm}
The Haantjes tensor associated with $V_Z$ vanishes for any choice
of the vector field $Z$.
\end{thm}

\noindent
\emph{Proof}. Let us write for simplicity $N$ in place of $N_{V_Z}$. Then, we have that
\begin{eqnarray*}
H_{V_Z}[X,Y]&=&N[Z\circ X,Z\circ Y]+Z^2\circ N[X,Y]-Z\circ N[X,Z\circ Y]-Z\circ N[Z\circ X,Y]=\\
&=&[X\circ Z^2,Z]\circ Y\circ Z-[X\circ Z,Z]\circ Z^2\circ Y+[Z,Y\circ Z^2]\circ X\circ Z+\\
&-&[Z,Y\circ Z]\circ X\circ Z^2+[X\circ Z]\circ Y\circ Z^2-[X,Z]\circ Z^3\circ Y+\\
&+&[Z,Y\circ Z]\circ X\circ Z^2-[Z,Y]\circ X\circ Z^3-[X\circ Z,Z]\circ Z^2\circ Y+\\
&+&[X,Z]\circ Z^3\circ Y-[Z,Y\circ Z^2]\circ X\circ Z+[Z,Y\circ Z]\circ X\circ Z^2+\\
&-&[X\circ Z^2,Z]\circ Y\circ Z+[X\circ Z,Z]\circ Z^2\circ Y-[Z,Y\circ Z]\circ X\circ Z^2+\\
&+&[Z,Y]\circ X\circ Z^3=0\ ,
\end{eqnarray*}
where $Z^2=Z\circ Z$ and $Z^3=Z\circ Z\circ Z$.
\nopagebreak
\unskip\nobreak\hskip2em plus1fill$\Box$\par\smallskip

Suppose now that $X$ is a vector field such that $V_X$ has everywhere distinct real eigenvalues $(v^1,\dots,v^n)$. Since the Haantjes tensor of $V_X$ vanishes, there exist local coordinates $(r^1,\dots,r^n)$ such that $(V_X)^i_j=\delta^i_jv^i$. These coordinates are the Riemann invariants of the corresponding system of hydrodynamic type. Moreover, we have

\begin{prop}
The components of the tensor field $c$ in the coordinates $(r^1,\dots,r^n)$ are given by
\[
c_{ij}^k=f_i\delta_i^k\delta_j^k\ .
\]
Furthermore, if $f_j\ne 0$ for all $j$,
then $f_i$ depends on the variable $r^i$ only.
\end{prop}

\noindent
\emph{Proof}.
In diagonal coordinates we have
$$\left(V_X\right)^i_j=c^i_{jk}X^k=v^i\delta^i_j,$$
hence, we get
$$c^j_{pq}c^i_{jk}X^k=c^j_{pq}v^i\delta^i_j=c^i_{pq}v^i.$$
On the other hand, due to the associativity of the algebra, we can also write
$$c^j_{pq}c^i_{jk}X^k=c^j_{pk}c^i_{jq}X^k=c^i_{jq}v^j\delta^j_p=c^i_{pq}v^p
\qquad\mbox{(no sum over $p$),}$$
and therefore,
\begin{equation*}
c^i_{pq}\left(v^i-v^p\right)=0.
\end{equation*}
Since the algebra is commutative and the eigenvalues of $V_X$ are pairwise distinct, this means that the structure constants, in the coordinates $(r^1,\dots,r^n)$, take the form
\beq\label{diagconpf}
c^i_{jk}=f_i\delta^i_j\delta^i_k,
\eeq
where the $f_i$ are arbitrary functions, depending in principle on all the variables $r^1,\dots,r^n$.
The requirement on the structure constants $c$ to satisfy the Herling-Manin condition \eqref{hmcomponents} implies further constraints on the functions $f_i$. Indeed, substituting \eqref{diagconpf} into \eqref{hmcomponents}, we get a set of equations the $f_i$ have to satisfy; considering for instance the case $m=j\neq k=i=l$, we get
$$f_j\partial_jf_k=0,$$
which means that, in the non-degenerate case when $f_j\neq 0$ for all $j$, then $f_k$ depends on $r^k$ only. It is easy to check that conditions \eqref{hmcomponents} give no further restrictions on the $f_i$; the proposition is proved.
\nopagebreak
\unskip\nobreak\hskip2em plus1fill$\Box$\par\smallskip

If the functions $f_i$ are everywhere different from zero, then it is easy to show that there exist local coordinates, called \emph{canonical coordinates}, such that $c_{ij}^k=\delta_i^k\delta_j^k$. Moreover, in this case, the vector field
$$
e=\sum_{i=1}^n\frac{1}{f_i}\frac{\partial}{\partial r^i}
$$
is globally defined and is the unity of the algebra.

\begin{rmk}
If the algebra has a unity $e$, then the Hertling-Manin condition implies
$${\rm Lie}_e c  
=0\ .$$
Indeed, for $X=Y=e$ the Hertling-Manin condition becomes
$$-[e,Z\circ W]+[e,Z]\circ W+[e,W]\circ Z=0.$$
\end{rmk}

\begin{rmk}
An alternative proof of the existence of canonical coordinates has been given in \cite{HM} under the assumption of semisimplicty of the algebra, that is, the existence of a basis of idempotents.
\end{rmk}

\section{Commutativity of the flows}
As a consequence of the Hertling-Manin condition, the conditions for the  commutativity of two hydrodynamical flows take a rather simple form.
\begin{prop}
\label{prop:comm}
The flows
\beq
\label{hdf1}
u^i_{t}=[V_X]^i_j u^j_x=c^i_{jk}X^j u^k_x
\eeq
and
\beq
\label{hdf2}
u^i_{\tau}=[V_Y]^i_j u^j_x=c^i_{jk}Y^j u^k_x
\eeq
commute if and only if the vector fields $X$ and $Y$ satisfy
 the condition
$$\left(\left({\rm Lie}_{X}c\right)(Y,Z)-\left({\rm Lie}_{Y}c\right)(X,Z)+[X,Y]\circ Z\right)\circ Z=0,$$
for any vector field $Z$. Equivalently,
\[
\begin{aligned}
\label{comcond-intri-polar}
& \left(\left({\rm Lie}_{X}c\right)(Y,Z)-\left({\rm Lie}_{Y}c\right)(X,Z)+[X,Y]\circ Z\right)\circ W\\
& \qquad
+\left(\left({\rm Lie}_{X}c\right)(Y,W)-\left({\rm Lie}_{Y}c\right)(X,W)+[X,Y]\circ W\right)\circ Z=0
\end{aligned}
\]
for all pairs $(Z,W)$ of vector fields.
In local coordinates this means that
$$
\begin{aligned}
& c^r_{is}\left[\left({\rm Lie}_{X}c\right)^i_{jq}Y^q-\left({\rm Lie}_{Y}c\right)^i_{jq}X^q+c^i_{jq}[X,Y]^q\right]\\
& \qquad +c^r_{ij}\left[\left({\rm Lie}_{X}c\right)^i_{sq}Y^q-\left({\rm Lie}_{Y}c\right)^i_{sq}X^q+c^i_{sq}[X,Y]^q\right]
=0\ .
\end{aligned}
$$
\end{prop}

\noindent
\emph{Proof}. It is well-known that the commutativity of the flows (\ref{hdf1}) and (\ref{hdf2}) is equivalent to the following requirements:

1. The $(1,1)$-tensor fields $V_{X}$ and
$V_{Y}$ (seen as endomorphism of the tangent bundle) commute.

2. For any vector field $Z$ the following condition is satisfied:
$$
[V_X(Z),V_Y(Z)]-V_X\left([Z,V_Y(Z)]\right)+V_Y\left([Z,V_X(Z)]\right)=0\ ,
$$
that is to say,
$$[Z\circ X,Z\circ Y]-X\circ[Z,Z\circ Y]+Y\circ[Z,Z\circ X]=0\ .$$
The first requirement is automatically verified due to the associativity of the algebra. Making use of identity (\ref{hmcond}), the second one
becomes
\beq
\label{commdot}
\left([Z\circ X,Y]+[X,Z\circ Y]-[X,Z]\circ Y-[X,Y]\circ Z-X\circ[Z,Y]\right)\circ Z=0.
\eeq
A simple calculation shows that the quantity in the bracket, namely
$$[Z\circ X,Y]+[X,Z\circ Y]-[X,Z]\circ Y-[X,Y]\circ Z-X\circ[Z,Y],$$
is equal to
\begin{equation}\label{qic}
\left({\rm Lie}_{X}c\right)(Y,Z)-\left({\rm Lie}_{Y}c\right)(X,Z)+[X,Y]\circ Z\ .
\end{equation}
Substituting \eqref{qic} into (\ref{commdot}), we get the thesis.
\nopagebreak
\unskip\nobreak\hskip2em plus1fill$\Box$\par\smallskip

\begin{cor}
A sufficient condition for the commutativity of the hydrodynamic flows (\ref{hdf1}) and (\ref{hdf2}) is
that
\begin{equation}\label{cc2-intri}
\left({\rm Lie}_{X}c\right)(Y,Z)-\left({\rm Lie}_{Y}c\right)(X,Z)+[X,Y]\circ Z=0
\end{equation}
for all vector fields $Z$, that is,
\beq\label{cc2}
\left({\rm Lie}_{X}c\right)^i_{pq}Y^q-\left({\rm Lie}_{Y}c\right)^i_{pq}X^q+c^i_{pq}[X,Y]^q=0
\eeq
or, equivalently,
\begin{equation}\label{cc3}
{\rm Lie}_{X}V_{Y}-{\rm Lie}_{Y}V_{X}-V_{[X,Y]}=0.
\end{equation}
\end{cor}

\section{Dubrovin principal hierarchy}
In this section, we adapt Dubrovin's construction of the principal hierarchy \cite{du93} to the case  of $F$-manifolds with compatible flat connection introduced by Manin in \cite{manin}.

\begin{defi}
An $F$-manifold with compatible flat connection is a manifold
 endowed with an associative commutative multiplicative structure given by a $(1,2)$-tensor field $c$ and 
 a flat torsionless connection $\nabla$ satisfying the symmetry condition
\beq\label{scc}
\nabla_l c^i_{jk}=\nabla_j c^i_{lk}\ ,
\eeq
meaning that $\nabla c$ is totally symmetric:
\beq\label{sccinv}
\left(\nabla_X c\right)\left(Y,Z\right)=\left(\nabla_Y c\right)\left(X,Z\right),
\eeq
for all vector fields $X$, $Y$, and $Z$.
\end{defi}
Notice that Hertling-Manin condition (\ref{hmcond}) does not appear in the above definition. Indeed,
 as proved by Hertling in \cite{hert}, it is a consequence of the existence of a torsionless (even non-flat) connection $\nabla$ satisfying (\ref{scc}). 

\begin{rmk}
 Notice that in flat coordinates condition \eqref{scc} reads
$$\d_l c^i_{jk}=\d_j c^i_{lk}.$$
This, together with the commutativity of the algebra, implies that
$$c^i_{jk}=\d_j C^i_k=\d_j\d_k C^i.$$
Therefore, condition (\ref{scc}) is equivalent to the local existence
 of a vector field $C$ satisfying, for any pair $(X,Y)$ of flat vector fields,
  the condition
$$X\circ Y=[X,[Y,C]].$$  
The above condition appears in the original definition of Manin \cite{manin}.
\end{rmk}

Let us construct now the principal hierarchy. In order to do so, the first step consists in defining the primary flows. Since the connection is flat, we can consider a basis $(X_{(1,0)},\dots, X_{(n,0)})$ of flat vector fields; the primary flows are thus defined as
\begin{equation}\label{primflo}
u^i_{t_{(p,0)}}=c^i_{jk}X^k_{(p,0)}u^j_x.
\end{equation}

\begin{prop}
The primary flows \eqref{primflo} commute.
\end{prop}

\noindent
\emph{Proof}.
Since the $X_{(p,0)}$ are flat and the torsion vanishes, they commute and
$$
{\rm Lie}_{X_{(p,0)}}c=\nabla_{X_{(p,0)}}c\ .
$$
Therefore, the commutativity condition \eqref{cc2-intri} for the vector fields $X=X_{(p,0)}$ and $Y=X_{(q,0)}$
follows from condition \eqref{scc}.
\nopagebreak
\unskip\nobreak\hskip2em plus1fill$\Box$\par\smallskip

Starting from the primary flows \eqref{primflo} one can introduce  the ``higher flows'' of the hierarchy, defined as
\beq
\label{hiflows}
u^i_{t_{(p,\alpha)}}=c^i_{jk}X^j_{(p,\alpha)}u^k_x,
\eeq
by means of the following recursive relations:
\beq\label{recrel}
\nabla_j X^i_{(p,\alpha)}=c^i_{jk}X^k_{(p,\alpha-1)}.
\eeq

\begin{rmk}
The flatness of the connection $\nabla$, the symmetry of the tensor $\nabla c$ (condition (\ref{scc}))
 and the associativity of the algebra with structure constants $c^i_{jk}$ are equivalent to the flatness of
the one-parameter family of connections
  defined, for any pair of vector fields $X$ and $Y$, by
$$\tilde{\nabla}_X Y=\nabla_X Y+z X\circ Y,\qquad z\in\mathbb{C}\ .$$
The vector fields obtained by means of the recursive relations (\ref{recrel}) are nothing but
the $z$-coefficients
 of a basis of flat vector fields of the deformed connection \cite{du93}.
\end{rmk}

In order to show that the higher flows (\ref{hiflows}) are well-defined, it is necessary to prove the following
\begin{prop}
\label{prop:rrcomp}
The recursive relations (\ref{recrel}) are compatible.
\end{prop}

\noindent
\emph{Proof}. We note that the recursive relations (\ref{recrel}) can be written in the form
$$\d_j X^i_{(p,\alpha)}=-\Gamma^i_{jk}X^k_{(p,\alpha)}-c^i_{kj}X^k_{(p,\alpha-1)},$$
thus, we have
\begin{eqnarray*}
(\d_j\d_m-\d_m\d_j)X^i_{(p,\alpha)}
&=&\left[\d_m\Gamma^i_{jl}-\d_j\Gamma^i_{ml}-
\Gamma^i_{jk}\Gamma^k_{ml}+\Gamma^i_{mk}\Gamma^k_{jl}\right]X^l_{(p,\alpha)}+\\
&&\left[\d_m c^i_{jl}-\d_j c^i_{ml}-
\Gamma^i_{kj}c^k_{ml}-\Gamma^k_{lm}c^i_{jk}+\Gamma^i_{km}c^k_{jl}+\Gamma^k_{lj}c^i_{mk}\right]X^l_{(p,\alpha-1)}\\
&&+\left[c^i_{jk}c^k_{ml}-c^i_{mk}c^k_{jl}\right]X^l_{(p,\alpha-2)}.
\end{eqnarray*}
The flatness of the connection $\nabla$, together with identity (\ref{scc}) and the associativity of the algebra, implies the vanishing of the quantity above. Therefore, relations \eqref{recrel} are compatible.
\nopagebreak
\unskip\nobreak\hskip2em plus1fill$\Box$\par\smallskip

Since the primary flows \eqref{primflo} commute and the recursive relations \eqref{recrel} are compatible, it only remains to prove the following
\begin{thm}
\label{thm:phcomm}
The flows of the principal hierarchy commute.
\end{thm}

\noindent
\emph{Proof}. Let us consider the hydrodynamic flows associated with the vector fields $X_{(p,\alpha)}$ and $X_{(q,\beta)}$. In order to show that these flows commute, we prove that they satisfy the sufficient condition \eqref{cc2}. In local coordinates it reads:
\begin{gather*}
X_{(p,\alpha)}^m(\d_m c^i_{jk})X_{(q,\beta)}^k
-X_{(q,\beta)}^m(\d_m c^i_{jk})X_{(p,\alpha)}^k+\\
-c^l_{jk}(\d_l X_{(p,\alpha)}^i)X_{(q,\beta)}^k
+c^i_{lk}(\d_j X_{(p,\alpha)}^l)X_{(q,\beta)}^k+\\
+c^i_{jl}(\d_k X_{(p,\alpha)}^l)X_{(q,\beta)}^k
+c^l_{jk}(\d_l X_{(q,\beta)}^i)X_{(p,\alpha)}^k+\\
-c^i_{lk}(\d_j X_{(q,\beta)}^l)X_{(p,\alpha)}^k
-c^i_{jl}(\d_k X_{(q,\beta)}^l)X_{(p,\alpha)}^k+\\
-c^i_{jk}\left((\d_l X_{(p,\alpha)}^k)X_{(q,\beta)}^l
+(\d_l X_{(q,\beta)}^k)X_{(p,\alpha)}^l\right)=0.
\end{gather*}
In particular, if the coordinates are flat, the first row vanishes due to the symmetry of the
tensor $\nabla c$. Moreover, using the recursive relations (\ref{recrel}) we obtain
\begin{gather*}
-c^l_{jk}c^i_{ln}X_{(p,\alpha-1)}^n X_{(q,\beta)}^k
+c^i_{lk}c^l_{jn} X_{(p,\alpha-1)}^n X_{(q,\beta)}^k+\\
+c^i_{jl}c^l_{kn} X_{(p,\alpha-1)}^n X_{(q,\beta)}^k
+c^l_{jk}c^i_{ln} X_{(q,\beta-1)}^n X_{(p,\alpha)}^k+\\
-c^i_{lk}c^l_{jn} X_{(q,\beta-1)}^n X_{(p,\alpha)}^k
-c^i_{jl}c^l_{kn} X_{(q,\beta-1)}^n X_{(p,\alpha)}^k+\\
-c^i_{jk}c^k_{mn} X_{(p,\alpha-1)}^n X_{(q,\beta)}^m
+c^i_{jk}c^k_{mn} X_{(q,\beta-1)}^n X_{(p,\alpha)}^m
\end{gather*}
which vanishes due to the associativity of the algebra.
\nopagebreak
\unskip\nobreak\hskip2em plus1fill$\Box$\par\smallskip

\begin{rmk}
The flows of the principal hierarchy are well-defined even in the case when
  the torsion of $\nabla$ does not vanish. However, their commutativity
   depends crucially on this additional assumption.
\end{rmk}

\section{$F$-manifolds with compatible connection and related integrable systems}

From the point of view of the theory of integrable systems of hydrodynamic type, the ``flat case'' and the associated principal hierarchy are exceptional. Therefore, it is quite natural to extend the notion of $F$-manifolds with compatible flat connection to the non-flat case. As a starting
  point, we consider an $F$-manifold endowed with a connection $\nabla$ satisfying
   (\ref{scc}). If $\nabla$ is flat, we know how to construct integrable systems of hydrodynamic type. Indeed, the starting point of the construction of the previous section
 is a basis of \emph{flat} vector fields, and the recursive procedure \eqref{recrel} defining the ``higher'' vector fields and the corresponding flows
  is well-defined as a consequence
 of the vanishing of the curvature. In the non-flat case, in order to define integrable systems of hydrodynamic type one needs to find an alternative
 way to select the vector fields.

\subsection{Hydrodynamic-type systems associated with $F$-manifolds}
In the flat case, the vector fields $X$ defining the principal hierarchy
 satisfy the condition
\beq\label{admvf-intri}
\left(\nabla_Z X\right)\circ W=\left(\nabla_W X\right)\circ Z
\eeq
for all pairs $(Z,W)$ of vector fields, that is, in local coordinates,
\beq\label{admvf}
c^i_{jm}\nabla_k X^m=c^i_{km}\nabla_j X^m\ .
\eeq
Indeed, in the case of the
 flat vector fields $X_{(p,0)}$ defining the primary flows, both sides of (\ref{admvf})
 vanish due to
$$\nabla_k X^m_{(p,0)}=0,\,\,\,p=1,\dots,n\ ,$$
while the vector fields defining the higher flows of the hierarchy satisfy \eqref{admvf} due to the associativity of the algebra:
$$c^i_{jm}\nabla_k X^m_{(p,\alpha)}=c^i_{jm}c^m_{kl}X^l_{(p,\alpha-1)}
=c^i_{km}c^m_{jl}X^l_{(p,\alpha-1)}=c^i_{km}\nabla_j X^m_{(p,\alpha)}.$$

A crucial remark is the following: if $\nabla$ satisfies condition (\ref{scc}),
 then \emph{any pair of solutions of (\ref{admvf}) defines commuting  flows even if the connection $\nabla$ is not flat}. More precisely, we have the following  
\begin{prop}
If $X$ and $Y$ are two vector fields satisfying condition
(\ref{admvf-intri}),
then the associated flows
\beq\label{flow1}
u^i_t=c^i_{jk}X^k u^j_x
\eeq
and
\beq\label{flow2}
u^i_{\tau}=c^i_{jk}Y^k u^j_x
\eeq
commute.
\end{prop}

\n
\emph{Proof}.
Recall from Proposition
\ref{prop:comm} that the flows (\ref{flow1})
 and (\ref{flow2}) commute if and only if
\beq
\label{comcond-intridimo}
\left(\left({\rm Lie}_{X}c\right)(Y,Z)-\left({\rm Lie}_{Y}c\right)(X,Z)+[X,Y]\circ Z\right)\circ Z=0
\eeq
for any vector field $Z$. On the other hand, the vanishing of the torsion of $\nabla$ gives the identity
$$
\left({\rm Lie}_{X}c\right)(Y,Z)=\left(\nabla_{X}c\right)(Y,Z)
-\nabla_{c(Y,Z)}X+c(Y,\nabla_Z X)+c(\nabla_Y X,Z)\ ,
$$
and this, together with the symmetry \eqref{sccinv} of $\nabla c$,  can be used to write the term in the bracket of
(\ref{comcond-intridimo}) as
$$
-\nabla_{Y\circ Z}X+\nabla_{X\circ Z}Y+[Y,X]\circ Z.
$$
Multiplying the above identity by $Z$, and using property (\ref{admvf-intri}) for the vector fields $X$ and $Y$, we obtain
\[
\begin{aligned}
&-\left(\nabla_{Y\circ Z}X\right)\circ Z+
\left(\nabla_{X\circ Z}Y\right)\circ Z+[Y,X]\circ Z^2=\\
&-\left(\nabla_{Z}X\right)\circ (Y\circ Z)+
\left(\nabla_{Z}Y\right)\circ (X\circ Z)+[Y,X]\circ Z^2=\\
&-\left(\nabla_{Y}X\right)\circ Z^2+
\left(\nabla_{X}Y\right)\circ Z^2+[Y,X]\circ Z^2=0.
\end{aligned}
\]
The proposition is proved.
\nopagebreak
\unskip\nobreak\hskip2em plus1fill$\Box$\par\smallskip

\begin{rmk}
From (\ref{scc}) and (\ref{admvf-intri}) it follows that the (1,1)-tensor field
$$(V_X)^i_j=c^i_{jk}X^k$$
satisfies the condition
$$\nabla_k (V_X)^i_j=\nabla_j (V_X)^i_k,$$
which is well-known in the Hamiltonian theory of systems of hydrodynamic type \cite{dn89}.
\end{rmk}

\subsection{Integrability condition}
In the flat case, we have seen that system (\ref{admvf}) admits a set of solutions, given by the vector fields of the principal hierarchy.  However, if $\nabla$ is non-flat, existence of solutions for system (\ref{admvf}) is not guaranteed; additional constraints have to be imposed on the curvature $R$ of the connection $\nabla$.

\begin{prop}
If $X$ is  a solution of (\ref{admvf-intri}), then the identity
\beq\label{sh}
Z\circ R(W,Y)(X)+
W\circ R(Y,Z)(X)+Y\circ R(Z,W)(X)=0,
\eeq
holds for any choice of the vector fields $(Y,W,Z)$.
\end{prop}

\n
\emph{Proof}. Condition (\ref{admvf-intri}) implies
$$\nabla_W(Z\circ\nabla_Y X-Y\circ\nabla_Z X)+\nabla_Y(W\circ\nabla_Z X-Z\circ\nabla_W X)+\nabla_Z(Y\circ\nabla_W X-W\circ\nabla_Y X)=0.$$
Using the symmetry condition (\ref{scc}) written in the form
$$\nabla_Y(X\circ Z)-\nabla_X(Y\circ Z)+Y\circ\nabla_X Z-X\circ\nabla_Y Z-[Y,X]\circ Z=0$$
we obtain identity (\ref{sh}).
\nopagebreak
\unskip\nobreak\hskip2em plus1fill$\Box$\par\smallskip

Condition (\ref{sh}) must be satisfied for \emph{any} solution $X$ of the
 system (\ref{admvf}). Since we are looking for a family of  
 vector fields satisfying (\ref{admvf}), it is natural to require that (\ref{sh}) holds true  for an \emph{arbitrary} vector field $X$. 
\begin{defi}
\label{defi:fmancc}
An \emph{$F$-manifold with compatible connection} is a manifold endowed with an associative commutative multiplicative structure given by a $(1,2)$-tensor field $c$
 and a torsionless connection $\nabla$ satisfying condition (\ref{sccinv})
 and condition
\beq\label{rc-intri}
Z\circ R(W,Y)(X)+
W\circ R(Y,Z)(X)+Y\circ R(Z,W)(X)=0,
\eeq
for any choice of the vector fields $(X,Y,W,Z)$.
In local coordinates this means that
\beq\label{shc}
R^k_{lmi}c^n_{pk}+R^k_{lip}c^n_{mk}+R^k_{lpm}c^n_{ik}=0\ .
\eeq
\end{defi}
 
\begin{rmk}
An equivalent form of condition (\ref{rc-intri}) can be easily obtained using the (second) Bianchi identity for the deformed connection
$$\tilde{\nabla}_X Y=\nabla_X Y+z X\circ Y,\qquad z\in\mathbb{C},$$
where $X$ and $Y$ are arbitrary vector fields. Indeed, due to associativity and symmetry condition (\ref{scc}), the Riemann tensor of this connection does not depend on $z$ \cite{St}.
 Using this fact it is easy to see that the Bianchi identity reduces to
\begin{eqnarray*}
0 &=& \tilde{\nabla}_X R(Y,Z)(W)+\tilde{\nabla}_Z R(X,Y)(W)+\tilde{\nabla}_Y R(Z,X)(W)\\
&=& X\circ R(Y,Z)(W)+Z\circ R(X,Y)(W)+Y\circ R(Z,X)(W)\\
&&\quad  - R(Y,Z)(X\circ W)- R(X,Y)(Z\circ W)-R(Z,X)(Y\circ W)
\end{eqnarray*}
for any choice of the vector fields $(X,Y,W,Z)$. Hence, condition \eqref{rc-intri} is equivalent to
$$R(Y,Z)(X\circ W)+R(X,Y)(Z\circ W)+R(Z,X)(Y\circ W)=0,$$
for every $(X,Y,W,Z)$.
\end{rmk}
From now on we will assume the existence of canonical coordinates $(r^1,\dots,r^n)$, discussing the meaning of condition (\ref{shc}) under this additional assumption. 

\begin{prop}
In canonical coordinates, system (\ref{admvf}) reduces to
\beq\label{chsym}
\d_k v^i=\Gamma^i_{ki}(v^k-v^i),\,\,\,\quad i\ne k\ ,
\eeq
where $v^i$ are the components of $X$ in such coordinates.
\end{prop}

\n
\emph{Proof}. Writing \eqref{admvf} in canonical coordinates, we get
$$\delta^i_j(\d_k v^i+\Gamma^i_{kl}v^l)=\delta^i_k(\d_j v^i+\Gamma^i_{jl}v^l).$$
In the case $i=j\ne k$, using the identities
\beq\label{egorov1}
\Gamma^i_{kk}=-\Gamma^i_{ki}
\eeq
and
\beq\label{egorov2}
\Gamma^i_{kl}=0,\qquad i\ne k\ne l\ne i,
\eeq
which follow from (\ref{scc}), we obtain system (\ref{chsym}). The remaining conditions give no further constraints.
\nopagebreak
\unskip\nobreak\hskip2em plus1fill$\Box$\par\smallskip

\begin{rmk}
We recall that, in canonical coordinates, the components of the vector field $X$  coincide with the characteristic velocities of the associated system of hydrodynamic type:
$$r^i_t=c^i_{jk}v^k r^j_x=v^i r^i_x\ ,\qquad i=1,\dots,n.$$
\end{rmk}
Compatibility conditions of system (\ref{chsym}) are well-known
 in the literature \cite{ts91}, and are given by the following conditions:
\begin{eqnarray}
\label{shprop1}
&&\d_i \Gamma^k_{mk}-\d_m\Gamma^k_{ik}=0,\\
\label{shprop2}
&&\d_i\Gamma^k_{km}-\Gamma^k_{km}\Gamma^m_{im}+\Gamma^k_{ik}\Gamma^k_{km}-\Gamma^k_{ik}\Gamma^i_{im}=0,
\end{eqnarray}
for pairwise distinct indices $k$, $i$, $m$.

\begin{prop}
Condition (\ref{shc}) is equivalent to conditions (\ref{shprop1}) and (\ref{shprop2}).
\end{prop}

\n
\emph{Proof}.
In canonical coordinates, condition (\ref{shc}) reads
\begin{eqnarray*}
&&R^k_{lmi}c^n_{pk}+R^k_{lip}c^n_{mk}+R^k_{lpm}c^n_{ik}=\\
&&R^k_{lmi}\delta^n_{p}\delta^n_{k}+R^k_{lip}\delta^n_{m}\delta^n_{k}+R^k_{lpm}\delta^n_{i}\delta^n_{k}=\\
&&R^n_{lmi}\delta^n_{p}+R^n_{lip}\delta^n_{m}+R^n_{lpm}\delta^n_{i}=0\ .
\end{eqnarray*}
If all the indices $m,i,p,n$ are distinct the above condition is trivially satisfied. Let us consider the case $n=p$ (the case $n\ne p$ can be treated in the same way
 and does not add further condition). If $n=i$, we obtain
$$R^n_{lmn}+R^n_{lnm}+\delta^n_m\,R^n_{lnn}=0\ ,$$
that is satisfied due to the skew-symmetry of the Riemann tensor with respect to the second and third lower indices.
 The same if $n=m$. For $n\ne i,m$, we obtain
\beq\label{rieone}
R^n_{nmi}=0,
\eeq
if $l=n$ and
\beq\label{riesec1}
R^n_{lmi}=0,
\eeq
if $l\ne n$. Since, due to (\ref{egorov1}), the components of the Riemann tensor
 vanish if all the indices are distinct, condition \eqref{riesec1} reduces to
\beq\label{riesec}
R^n_{mmi}=0,\qquad n\ne m\ne i\ne n.
\eeq
Finally, using (\ref{egorov1}) and  (\ref{egorov2}), it is easy to check that conditions (\ref{rieone}) and (\ref{riesec}) 
 are equivalent to conditions (\ref{shprop1}) and (\ref{shprop2}) respectively. This proves the proposition.
\nopagebreak
\unskip\nobreak\hskip2em plus1fill$\Box$\par\smallskip

\begin{rmk}
If the compatibility conditions (\ref{shprop1}) and (\ref{shprop2}) are satisfied, the general solution of the system (\ref{chsym})
 depends on $n$ arbitrary functions of a single variable. Moreover, due to (\ref{shprop1}),  
 any solution $(v^1,\dots,v^n)$ of
 (\ref{chsym}) satisfies the condition
\beq\label{shprop}
\d_k\left(\frac{\d_j v^i}{v^j-v^i}\right)-\d_j\left(\frac{\d_k v^i}{v^k-v^i}\right)=0,\,\,\,i\ne j\ne k\ne i,
\eeq
known in literature as semi-Hamiltonian property \cite{ts91}. An invariant and highly non trivial formulation of such a property has been found in \cite{PSS}.
\end{rmk}

Due to the above remark, under the assumption of existence of canonical coordinates we have a set of solutions of (\ref{chsym}) leading to a family of commuting systems of hydrodynamic type, depending on $n$ arbitrary functions. This result shows the deep relation between $F$-manifold with compatible connection (Definition \ref{defi:fmancc}) and integrable systems of PDEs.

\section{Riemannian $F$-manifolds and Egorov metrics}\label{egorovsection}

In this section we consider the special case where the connection $\nabla$ is a metric connection. This assumption plays an important role in the Hamiltonian theory of systems of hydrodynamic type (see for instance \cite{du90,Pa2007,PT} and references therein), as well as in the theory of Frobenius manifolds \cite{du93,du97}.

\begin{defi}
A \emph{Riemannian $F$-manifold} is an $F$-manifold with a compatible connection $\nabla$ satisfying the following additional conditions:

\n
1. The connection is metric: $$\nabla g=0\ .$$

\n
2. The inner product $\langle\cdot,\cdot\rangle$ defined by the metric $g$ is invariant with respect to the product $\circ$:
\beq\label{invpr}
\langle X\circ Y,Z\rangle=\langle X,Y\circ Z\rangle\ .
\eeq
\end{defi}
In local coordinates, condition (\ref{invpr}) reads
\beq\label{invpr2}
g_{iq}c^q_{lp}=g_{lq}c^q_{ip},\quad \mbox{or}\quad g^{iq}c^l_{qp}=g^{lq}c^i_{qp},
\eeq
where $g_{ij}$ and $g^{ij}$ are respectively  the covariant and the contravariant components of the metric $g$.

If there exist canonical coordinates, the metric $g$ entering the definition
 of Riemannian $F$-manifold is an Egorov metric. Let us recall the definition of this special class of metrics.

\begin{defi}
A metric is called \emph{Egorov} if there exist coordinates $(r^1,\dots,r^n)$ such that it is diagonal and \emph{potential}:
$$g_{ij}=\delta^i_j\,g_{ii}(r^1,\dots,r^n)=\delta^i_j\,\d_i F,$$
for a certain function $F$.
\end{defi}

Now, if we assume the existence of canonical coordinates, condition (\ref{invpr2}) tells us that the metric $g$ is diagonal in such coordinates,
 while condition (\ref{egorov2})---which follows from (\ref{scc})---implies that the metric is potential. Therefore, $g$ is an Egorov metric. Conversely,
 given an Egorov metric $g$ whose curvature tensor satisfies condition \eqref{riesec}, we can locally construct a 
 Riemannian $F$-manifold.
 More precisely, let  $\left(r^1,\dots,r^n\right)$ be the coordinates where  $g$ is diagonal and potential.
  Then, the metric $g$ and the structure constants
$$c^i_{jk}(r)=\delta^i_j\delta^i_k$$
endow the open set where the coordinates $\left(r^1,\dots,r^n\right)$ are defined with the structure of a Riemannian $F$-manifold. 

We point out that condition \eqref{shc} is far from being trivial. Indeed, using the above remark, it is easy to construct examples of metrics satisfying  properties (\ref{invpr}) and (\ref{scc}). Much more difficult is the problem of finding Egorov metrics
 which satisfy also condition (\ref{shc}), since the potential has to fulfill (\ref{riesec}).
 However, there exists an important class of metrics, appearing in the
Hamiltonian theory of integrable hierarchies of hydrodynamic type (not necessarily of Egorov type) whose curvature satisfies \eqref{shc}. These are the metrics whose Riemann tensor admits ``a quadratic expansion'' in terms of the flows of the hierarchy \cite{fe91,fm90}:
$$u^i_{t_{\alpha}}=c^i_{jk}X^k_{(\alpha)}u^j_x,\,\,\,\quad i=1,\dots,n.$$
This means that
\beq\label{qexp}
R^{sk}_{mi}=\left(c^s_{ml}c^k_{iq}-c^s_{il}c^k_{mq}\right)\sum_{\alpha}\epsilon_{\alpha}X^l_{(\alpha)}X^q_{(\alpha)},\qquad\epsilon_{\alpha}=\pm 1,
\eeq
where the index $\alpha$ can take value on a finite or infinite---even continuous---set. 

\begin{prop}\label{quadexpprop}
Suppose that $\nabla$ is the Levi-Civita connection of a metric $g$,  and that its curvature satisfies condition (\ref{qexp}).
 In this case, condition (\ref{shc}) is automatically satisfied.
\end{prop}

\n
\emph{Proof}. We have that
\begin{eqnarray*}
&&R^{sk}_{mi}c^n_{pk}+R^{sk}_{ip}c^n_{mk}+R^{sk}_{pm}c^n_{ik}=\\
&&\sum_{\alpha}\epsilon_{\alpha}[(c^s_{mr}c^k_{iq}-c^s_{ir}c^k_{mq})c^n_{pk}+(c^s_{ir}c^k_{pq}-c^s_{pr}c^k_{iq})c^n_{mk}
+(c^s_{pr}c^k_{mq}-c^s_{mr}c^k_{pq})c^n_{ik}]X^r_{(\alpha)}X^q_{(\alpha)}=\\
&&\sum_{\alpha}\epsilon_{\alpha}[(c^k_{iq}c^n_{pk}-c^k_{pq}c^n_{ik})c^s_{mr}+(c^k_{pq}c^n_{mk}-c^k_{mq}c^n_{pk})c^s_{ir}
+(c^k_{mq}c^n_{ik}-c^k_{iq}c^n_{mk})c^s_{pr}]X^r_{(\alpha)}X^q_{(\alpha)}\ ,
\end{eqnarray*}
which vanishes due to associativity.
\nopagebreak
\unskip\nobreak\hskip2em plus1fill$\Box$\par\smallskip


\begin{rmk}
If the functions
$$g^{lq}:=\sum_{\alpha}\epsilon_{\alpha}X^l_{(\alpha)}X^q_{(\alpha)}$$
define the contravariant components of a metric satisfying condition (\ref{invpr2}), then the operator
$$\sum_{\alpha}\epsilon_{\alpha}\left(w_{\alpha}\right)^i_ku^k_x
\left(\frac{d}{dx}\right)^{\!-1}\!\!\!\left(w_{\alpha}\right)^j_hu^h_x,\qquad\left(w_{\alpha}\right)^i_j:=c^i_{jk}X^k_{(\alpha)}$$
is a purely nonlocal Poisson operator (see \cite{GLR2} for details).
\end{rmk}

\section{An example: reductions of the dispersionless KP hierarchy}
In this section we will consider a class of Riemannian $F$-manifolds associated with a well-known class of hydrodynamic type systems: the reductions of the dispersionless KP hierarchy. For a generic reduction, the metric will be non-flat.

The dispersionless KP (or dKP) hierarchy can be defined by introducing the formal series
\begin{equation}\label{hydronor}
\lambda=p+\sum_{k=0}^\infty\frac{A^k}{p^{k+1}},
\end{equation}
which has to satisfy the following dispersionless Lax equations
\begin{equation*}
\lambda_{t_n}=\left\{\lambda,\, \frac{1}{n}\left(\lambda^n\right)_+\right\}.
\end{equation*}
Here $\{f,g\}=\partial_xf\,\partial_pg-\partial_pf\,\partial_xg$ denotes the canonical Poisson bracket, and $(\,\cdot\,)_+$ is the polynomial part of the argument. For simplicity, we will consider here only the second flow ($n=2$); all other flows of the hierarchy can be treated in the same way. For the second flow, we have
\begin{equation}\label{laxbenney}
\lambda_{t_2}=\left\{\lambda,\, \frac{1}{2}p^2+A^0\right\}=p\lm_x-A^0_x\lm_p,
\end{equation}
or, explicitly in terms of the variables $A^k$,
\beq
\label{benneychain}
A^k_{t_2}=A^{k+1}_x+kA^{k-1} A^0_x, \qquad k=0,1,2,\dots
\eeq
This last system is also known in the literature as Benney chain \cite{be73}; its Lax representation \eqref{laxbenney} appeared for the first time in \cite{lema}.   An $n$-component reduction of the Benney chain is a restriction of the infinite dimensional system (\ref{benneychain}) to a suitable $n$-dimensional submanifold, that is
$$A^k = A^k(u^1,\dots,u^n),\qquad k = 0, 1,\dots\,.$$
The reduced systems are systems of hydrodynamic type in the variables $(u^1,\dots, u^n)$ that parametrize the submanifold:
$$u^i_t = v^i_j(u)u^j_x,\qquad i = 1,\dots,n.$$
Reductions of the Benney system were introduced in \cite{GT}, and there it was proved that
such systems are diagonalizable and integrable via the generalized hodograph transformation \cite{ts91}. Clearly, in the case of a reduction, the coefficients of the series (\ref{hydronor}) depend on the Riemann invariants $(r^1,\dots , r^n)$
  and the series can be thought as the asymptotic
expansion for $p\to\infty$ of a suitable function $\lm(p,r^1,\dots,r^n)$ depending piecewise analytically on the parameter $p$. It turns out \cite{GT,GT2} that such a function satisfies a system of chordal Loewner equations,
\begin{equation}\label{loewner}
\frac{\de\lambda}{\de r^i}=\frac{\de_i A^0}{p-v^i}\,\lambda_p, \qquad i=1\dots,n,
\end{equation}
describing families of conformal maps (with respect to $p$) in the complex upper half plane. The analytic properties of $\lm$ characterize the
reduction. More precisely, in the case of an $n$-reduction the associated function $\lm$ possesses $n$ distinct critical points on the real axis; these are the characteristic velocities $v^i$ of the reduced system, that is,
$$\lm_p(v^i):=\frac{\de\lm}{\de p}(v^i)=0, \qquad i=1,\dots,n,$$
and the corresponding critical values can be chosen as Riemann invariants. Compatibility conditions of the Loewner system \eqref{loewner} are of the form
\begin{align*}
\de_iv^j&=\frac{\de_iA^0}{v^i-v^j}\notag\\
&\hspace{3cm} i\neq j,\\
\de^2_{ij}A^0&=\frac{2\de_iA^0\de_jA^0}{(v^i-v^j)^2}\notag
\end{align*}
and were found by Gibbons and Tsarev \cite{GT2}. Thus, every reduction of the Benney chain is described by a particular solution of the Loewner system \eqref{loewner}.

Starting from the function $\lm$, we will show now how to give to the manifold parametrized by the Riemann invariants $\left(r^1,\dots,r^n\right)$, a structure of $F$-manifold with a compatible connection---in general non-flat. In order to do this, we define a metric
\begin{equation}\label{metricbenney}
g(\partial,\partial')=\sum_{i=1}^n\underset{p=v^i}{\rm{res}}
\left(\frac{\partial\lambda(p)\,\partial'\lambda(p)}{\lambda_p}dp\right),
\end{equation}
and structure constants
\begin{equation}\label{constantsbenney}
c(\partial,\partial',\partial'')=\sum_{i=1}^n\underset{p=v^i}{\rm{res}}
\left(\frac{\partial\lambda(p)\,\partial'\lambda(p)\,
\partial''\lambda(p)}{\lambda_p}dp\right),
\end{equation}
where $\partial$, $\partial'$, $\partial''$ are arbitrary tangent vectors on the manifold. In the coordinates $\left(r^1,\dots,r^n\right)$, and making use of the Loewner equations \eqref{loewner}, the metric takes the diagonal form
\begin{align*}
g\left(\frac{\partial}{\partial r^i},\frac{\partial}{\partial r^j}\right)&=\sum_{i=1}^n\underset{p=v^i}{\rm{res}}
\left(\frac{\partial\lambda}{\partial r^i}\frac{\partial\lambda}{\partial r^j}\frac{dp}{\lambda_p}\right)=\sum_{i=1}^n\underset{p=v^i}{\rm{res}}
\left(\frac{\de_i A^0\de_j A^0\,\,\lambda_p\, dp}{(p-v^i)(p-v^j)}\right)\\
&\\
&=\de_i A^0\de_j A^0\lambda_{pp}(v^i)\,\delta_{ij}=\de_i A^0\,\delta_{ij},
\end{align*}
where we used the fact \cite{glr2008} that
$$\lambda_{pp}(v^i)=\frac{1}{\partial_i A^0}.$$
In particular, the metric is Egorov. Moreover, a similar calculation for the structure constants gives
\begin{align*}
c\left(\frac{\partial}{\partial r^i},\frac{\partial}{\partial r^j},\frac{\partial}{\partial r^k}\right)&=\sum_{i=1}^n\underset{p=v^i}{\rm{res}}
\left(\frac{\partial\lambda}{\partial r^i}\frac{\partial\lambda}{\partial r^j}\frac{\partial\lambda}{\partial r^k}\frac{dp}{\lambda_p}\right)=\sum_{i=1}^n\underset{p=v^i}{\rm{res}}
\left(\frac{\de_i A^0\de_j A^0\de_k A^0\,\,\left(\lambda_p\right)^2\,dp}{(p-v^i)(p-v^j)(p-v^k)}\right)\\
&\\
&=\de_i A^0\de_j A^0\de_k A^0\,\,\left(\lambda_{pp}(v^i)\right)^2\,\delta_{ij}\delta_{ik}=\de_i A^0\,\delta_{ij}\delta_{ik},
\end{align*}
and from this it follows that
$$\frac{\partial}{\partial r^i}\circ\frac{\partial}{\partial r^j}=\delta_{ij}\frac{\partial}{\partial r^i},$$
namely $\left(r^1,\dots,r^n\right)$ are canonical coordinates for the algebra.
\begin{rmk}
The metric \eqref{metricbenney} and the structure constants \eqref{constantsbenney} were introduced for the first time by Dubrovin in \cite{du93}, in the particular case of the Gelfand-Dikii reductions of the dKP hierarchy, where the function $\lambda$ is a polynomial in $p$. The same metric and constants were also used by    Chang \cite{Chang2007} and Ferguson and Strachan \cite{FS}, for the study of reductions where $\lambda$ is rational or logarithmic. We remark that in all these examples the metric considered turns out to be flat.  
\end{rmk}
We have now to prove that the metric and the structure constants defined in this way are compatible, namely that conditions \eqref{scc} and \eqref{shc} are satisfied. As regard condition \eqref{scc}---due to the results of Section \ref{egorovsection}---it is sufficient to note that the metric \eqref{metricbenney} is Egorov. On the other hand, for condition \eqref{shc},  we only have to recall the result of \cite{glr2008}, where the curvature tensor of the metric \eqref{metricbenney} has been shown to possess the following quadratic expansion:
$$R^{ij}_{ij}=\frac{1}{2\pi i}\int_{C}w^i(\lm)\,w^j(\lm)\,d\lm,\qquad w^i(\lambda)=\frac{\frac{\de
p}{\de\lm}}{(p(\lm)-v^i)^2},$$
where $p(\lm)=\lm^{-1}(p)$ is the inverse of $\lm$ with respect to $p$, and $C$ is a suitable contour on the complex $\lm$-plane.
Due to Proposition \ref{quadexpprop}, the existence of a quadratic expansion of the curvature implies that condition \eqref{shc} is satisfied. Alternatively, such a condition follows from the well-known fact that the characteristic velocities $v^i$---which satisfy condition \eqref{chsym}---satisfy the semi-Hamiltonian condition \eqref{shprop}.

\begin{rmk}
A similar construction can be done using instead of the metric (\ref{metricbenney}),
 one of the metrics
\begin{equation}\label{metricbenney2}
g\left(\frac{\partial}{\partial r^i},\frac{\partial}{\partial r^j}\right)=\sum_{i=1}^n\underset{p=v^i}{\rm{res}}
\varphi_i(r^i)\left(\frac{\partial\lambda}{\partial r^i}\frac{\partial\lambda}{\partial r^j}\frac{dp}{\lambda_p}\right)
\end{equation}
where $\varphi_i$ are arbitrary functions of a single variable, and defining the corresponding structure constants as
\begin{equation}\label{constantsbenney2}
c\left(\frac{\partial}{\partial r^i},\frac{\partial}{\partial r^j},\frac{\partial}{\partial r^k}\right)=\sum_{i=1}^n\underset{p=v^i}{\rm{res}}
(\varphi_i(r^i))^2\left(\frac{\partial\lambda}{\partial r^i}\frac{\partial\lambda}{\partial r^j}\frac{\partial\lambda}{\partial r^k}\frac{dp}{\lambda_p}\right).
\end{equation}
If all the functions $\varphi_i$ are different from zero, it turns out that the structure constants (\ref{constantsbenney2}) admit canonical coordinates.
 Moreover, in such coordinates, the metric (\ref{metricbenney2}) is potential.
 In this way, repeating the construction described in this section, one defines, for any choice of the functions $\varphi_i$,
  a new structure of $F$-manifold with compatible connection on the \emph{same} manifold.  Notice that in case of  reductions related to Frobenius manifolds, such as the Zakharov and the Gel'fand-Dikii reductions \cite{du93,dlz08}, one of the metrics (\ref{metricbenney2}) is the \emph{intersection form} of the Frobenius manifold. Using this metric, the construction above reduces to the \emph{Dubrovin's duality} of the theory of Frobenius manifolds \cite{du04}.
\end{rmk}

\subsection*{Acknowledgments}
We thank Maxim Pavlov for stimulating discussions. M.P.\ and A.R.\ would like to thank the Department {\em Matematica e
Applicazioni\/} of the Milano-Bicocca University for the hospitality. Support of A.R. through the ESF grant MISGAM 2265 for his visit to the Milano-Bicocca department is gratefully acknowledged.


\begin{thebibliography}{99}
\bibitem{be73} D.J. Benney,
\emph{Some properties of long nonlinear waves},
Stud. Appl. Math. {\bf 52} (1973), 45--50.

\bibitem{Chang2007}
Jen-Hsu Chang, \emph{On the waterbag model of the dispersionless KP hierarchy (II)}, J. Phys. A: Math. Theor. {\bf 40} (2007), 12973--12985.

\bibitem{du90}
B.A. Dubrovin, \emph{On the differential geometry of strongly integrable systems of hydrodynamics type}, (Russian)  Funktsional. Anal. i Prilozhen. {\bf 24} (1990),  no. 4, 25--30, 96;  translation in  Funct. Anal. Appl. {\bf 24}  (1990),  no. 4, 280--285 (1991).

\bibitem{du93}
B.A. Dubrovin, \emph{Geometry of 2D topological field theories},
in: Integrable Systems and Quantum Groups, Montecatini Terme, 1993.
Editors: M. Francaviglia, S. Greco. Springer Lecture Notes in
Math. {\bf 1620} (1996), pp.\ 120--348.


\bibitem{du97}
B.A. Dubrovin, \emph{Flat pencils of metrics and {F}robenius manifolds},
in: Integrable systems and algebraic geometry ({K}obe/{K}yoto,1997)
World Sci. Publ., River Edge, NJ {\bf 1620} (1998), pp.\ 47--72.


\bibitem{du04}
B.A. Dubrovin, \emph{On almost duality for Frobenius manifolds}, in: Geometry, topology, and mathematical physics,  Amer. Math. Soc. Transl. Ser. 2, 212, Amer. Math. Soc., Providence, RI, 2004, pp. 75--132.


\bibitem{dlz08}
B.A. Dubrovin, S.Q. Liu, Y. Zhang, \emph{Frobenius manifolds and central invariants for the {D}rinfeld-{S}okolov bi{H}amiltonian structures}, Adv. Math. {\bf 219} (2008),  no. 3, 780--837.


\bibitem{dn89} B.A. Dubrovin, S.P. Novikov, \emph{Hydrodynamics of weakly deformed
soliton lattices. Differential geometry and Hamiltonian theory}, Uspekhi Mat.
Nauk {\bf 44} (1989), 29--98. English translation in Russ. Math. Surveys {\bf 44}
(1989), 35--124.

\bibitem{fe91} E.V. Ferapontov,
\emph{Differential geometry of nonlocal Hamiltonian operators
 of hydrodynamic type}, Funct. Anal. Appl. {\bf 25}  (1991),
  no. 3, 195--204.

\bibitem{fm90} E.V. Ferapontov, O.I. Mokhov,
\emph{Nonlocal Hamiltonian operators of hydrodynamic type that
 are connected with metrics of constant curvature},
  Russ. Math. Surv. {\bf 45} (1990),  no. 3, 218--219.

\bibitem{FS} J.T. Ferguson, I.A.B. Strachan, \emph{Logarithmic deformations of the rational superpotential/Landau-Ginzburg construction of solutions of the WDVV equations}, Comm.\ Math.\ Phys.\  
{\bf 280} (2008), 1--25.

\bibitem{glr2008} J. Gibbons, P. Lorenzoni, A. Raimondo,
\emph{Hamiltonian structure of reductions of the Benney system},
Comm.\ Math.\ Phys.\ {\bf 287} (2009), 291-322,

\bibitem{GLR2} J. Gibbons, P. Lorenzoni, A. Raimondo,
\emph{Purely nonlocal Hamiltonian formalism for systems of hydrodynamic type},
 arXiv:0812.3317.

\bibitem{GT} J. Gibbons, S.P. Tsarev,
\emph{Reductions of the Benney equations},
 Phys. Lett. A {\bf 211} (1996), no. 1, 19--24.
 
\bibitem{GT2}  J. Gibbons, S.P. Tsarev,
\emph{Conformal maps and reductions of the Benney equations},
  Phys. Lett. A  {\bf 258} (1999),  no. 4-6, 263--271.

\bibitem{hert}
C. Hertling, \emph{Multiplication on the tangent bundle},
 arXiv:math/9910116 
 
\bibitem{HM}
C. Hertling, Y. Manin, \emph{Weak Frobenius manifolds}, Internat. Math. Res. Notices {\bf 1999}, no. 6, 277--286.

\bibitem{KM}
B.G. Konopelchenko, F. Magri, \emph{Coisotropic deformations of associative algebras and dispersionless integrable hierarchies}, Comm. Math. Phys. {\bf 274} (2007), 627--658.

\bibitem{lema}
D. Lebedev, Y. Manin, \emph{Conservation laws and Lax representation of Benney's long wave equations},
Phys. Lett. A {\bf 74} (1979), 154--156.

\bibitem{manin} Y. Manin,
\emph{$F$-manifolds with flat structure and Dubrovin's duality}, Adv. Math. {\bf 198} (2005), no. 1, 5--26.

\bibitem{Pa2007}
M.V. Pavlov, \emph{Integrability of Egorov systems of hydrodynamic type}, (Russian)  Teoret. Mat. Fiz.
{\bf 150} (2007),  no. 2, 263--285;  translation in  Theoret. and Math. Phys. {\bf 150} (2007),  no. 2, 225--243.

\bibitem{PSS} M.V. Pavlov, S.I. Svinolupov, R.A. Sharipov, \emph{An invariant criterion for hydrodynamic integrability}, (Russian)  Funktsional. Anal. i Prilozhen. {\bf 30} (1996),  no. 1, 18--29, 96;  translation in  Funct. Anal. Appl. {\bf 30} (1996),  no. 1, 15--22.

\bibitem{PT} M.V. Pavlov, S.P. Tsarev,  \emph{Tri-Hamiltonian structures of Egorov systems of hydrodynamic type}, (Russian)  Funktsional. Anal. i Prilozhen. {\bf 37} (2003),  no. 1, 38--54.

\bibitem{St}
I.A.B Strachan, \emph{Frobenius manifolds: natural submanifolds and induced bi-Hamiltonian structures},  Differential Geom. Appl. {\bf 20} (2004),  no. 1, 67--99.

\bibitem{ts91} S.P. Tsarev,
\emph{The geometry of Hamiltonian systems of hydrodynamic type. The
generalised hodograph transform},
USSR Izv. {\bf 37} (1991) 397--419.

\bibitem{za80} V.E. Zakharov,
\emph{Benney equations and quasiclassical approximation in the inverse problem},
 Funktional. Anal. i Prilozhen {\bf 14} (1980), 15--24.

\end{thebibliography}
\end{document}